\documentclass[11pt]{article}
\usepackage{amsmath}
\usepackage{amssymb}

{\catcode `\@=11 \global\let\AddToReset=\@addtoreset}
\AddToReset{equation}{section}

\usepackage{epsfig}
\usepackage{amsmath}
\usepackage{amssymb}
\usepackage{graphicx}
\usepackage{color}

\usepackage[latin1]{inputenc}

\newtheorem{theorem}{Theorem}[section]
\newtheorem{lemma}{\bf Lemma}[section]

\newtheorem{proposition}{Proposition}[section]
\newtheorem{@definition}{\sc Definition}[section]

\newtheorem{@remark}{\sc Remark}[section]

\newtheorem{@example}{\sc Example}[section]

\newcommand{\beqn}{\begin{displaymath}}
\newcommand{\eeqn}{\end{displaymath}}
\newcommand{\beq}{\begin{equation}}  
\newcommand{\eeq}{\end{equation}}

\def\mathsf{\bf}
\def\N{\mathbb{N}}

\def\R{\mathbb{R}}
\def\Z{\mathbb{Z}}

\def\i{\mathrm i}
\def\d{\mathrm d}
\def\e{\mathrm e}

\def\cH{\cal{H}}

\def\E{\mathrm E}

\def\text{\mbox}

\def\1{{\bf 1}}

\newcommand{\mbf}[1]{\mbox{\boldmath $#1$}}

\newcommand{\nn}{\nonumber}
\newcommand{\noi}{\noindent}

\newcommand{\mbx}{{\mbf x}}

\newcommand{\mbgamma}{{\mbf \gamma}}

\def\eq2{
\stackrel{\small \rm mod \,2}{=}}

\def\n2{
\stackrel{\small \rm mod \,2}{\neq}}

\def\limfdd{\renewcommand{\arraystretch}{0.5}
\begin{array}[t]{c}
\stackrel{\rm fdd}{\longrightarrow} \\
\end{array}\renewcommand{\arraystretch}{1}}

\def\eqfdd{\renewcommand{\arraystretch}{0.5}
\begin{array}[t]{c}
\stackrel{\rm fdd}{=} \\
\end{array}\renewcommand{\arraystretch}{1}}

\def\neqfdd{\renewcommand{\arraystretch}{0.5}
\begin{array}[t]{c}
\stackrel{\rm fdd}{\neq} \\
\end{array}\renewcommand{\arraystretch}{1}}

\newtheorem{thm}{Theorem}[section]

\newtheorem{defn}[thm]{Definition}
\newtheorem{rem}{Remark}[section]

\newtheorem{ex}[thm]{Example}

\def\vep{\varepsilon}

\begin{document}

\title{Scaling  transition
for long-range dependent \\
Gaussian random fields}


\author{Donata Puplinskait\.e
and Donatas Surgailis\footnote{Corresponding author. E-mail: donatas.surgailis@mii.vu.lt}
\\
\small \it Vilnius University}
\maketitle

\begin{abstract}
In \cite{ps2014} we introduced the notion of scaling  transition
for stationary random fields $ X  $
on $\Z^2$
in terms of partial sums limits, or scaling limits,
of $X$ over rectangles
whose sides grow at possibly different rate.
The present paper establishes the existence of scaling  transition for a natural  class
of stationary Gaussian random fields  on $\Z^2$ with long-range dependence.
The scaling limits
of such random fields are identified
and characterized by dependence properties of rectangular increments.

\end{abstract}

\smallskip
{\small

\noi {\it Keywords:} scaling  transition; long-range dependence;
Gaussian random field;
operator scaling random field

}

\vskip.7cm

\section{Introduction}

Let $X = \{ X(t,s); (t,s) \in \Z^2\}$ be a stationary random field (RF) on the lattice $\Z^2$, $\gamma >0$ a given number
and
$K_{[nx, n^{\gamma}y]} := \{ (t,s)  \in \Z^2: 1\le t \le nx, 1\le s \le n^{\gamma}y \} $ be a sequence of rectangles whose sides grow at possibly different
rate  $O(n)$ and $O(n^{\gamma})$. Assume that for any $\gamma >0$ there exist
a nontrivial RF $V_\gamma =
\{V_\gamma (x,y); (x,y) \in \R^2_+ \}$ and a normalization $A_n(\gamma)  \to \infty $
such that
\begin{equation} \label{Xsum01}
A^{-1}_n(\gamma) \sum_{(t,s) \in K_{[nx, n^{\gamma}y]}} X(t,s) \ \limfdd \ V_\gamma (x,y), \quad (x,y) \in \R^2_+, \quad n \to \infty.
\end{equation}
We say that RF {\it $X$ exhibits scaling  transition} if there exists $\gamma_0 >0$ such that
\begin{equation}\label{Vtrans}
V_\gamma \eqfdd V_+, \ \gamma >\gamma_0, \quad V_\gamma \eqfdd V_-, \ \gamma <\gamma_0 \quad \text{and} \quad
V_+ \neqfdd a V_- \ (\forall \, a >0).
\end{equation}
See the end of this sec. for all unexplained notation.
In other words, \eqref{Vtrans} say that the
scaling limits $V_\gamma $ in \eqref{Xsum01}
do not depend on $\gamma $ for  $\gamma >\gamma_0$ and $\gamma < \gamma_0$ and are different up to a multiplicative constant
(the last condition is needed to exclude a trivial change of the scaling limit by a linear  change of normalization).
In the sequel, $V_{\gamma_0}$ will be called the {\it well-balanced} scaling limit of $X$, and
$V_+, V_- $ the {\it unbalanced} scaling limits of $X$.
Obviously,
if the limits $V_\gamma \eqfdd V $ in \eqref{Xsum01}
are  the same for any $\gamma >0$, the RF $X$ does not exhibit  scaling  transition.

The notion of scaling  transition was introduced in our paper \cite{ps2014}, which also established
the existence of such transition for a class of aggregated $\alpha$-stable $(1< \alpha \le 2)$ RFs on $\Z^2$. The last paper
identified the scaling limits $V_+, V_-, V_{\gamma_0}$ and characterized these RFs  by certain dependence properties
of increments on rectangles $K \subset \R^2_+$.

The present paper extends the results of \cite{ps2014} by proving the existence of  scaling  transition
for a natural class of stationary long-range dependent (LRD) Gaussian RFs with Type I spectral density $f_{\rm I}$
and its absence for  Gaussian RFs with Type II spectral density $f_{\rm II}$ in \eqref{typeI/IIsp}:
\begin{equation} \label{typeI/IIsp}
\text{Type I density:}\quad f_{\rm I}(x,y)= \frac{g(x,y)}{\big(|x|^2 + c|y|^{2 H_2/H_1}\big)^{H_1/2}},
\qquad \text{Type II density:} \quad  f_{\rm II}(x,y)= \frac{g(x,y)}{|x|^{2d_1}|y|^{2 d_2}},
\end{equation}
where $H_1, H_2>0, H_1 H_2 < H_1+H_2,  \, c>0, \, 0<  d_1, d_2 < 1/2$ are parameters
and $g$ is a bounded positive function having a positive limit $g(0,0)>0$  at the origin (w.l.g., we assume $g(0,0) =1 $).
Type II spectral densities $ f_{\rm II}$ in \eqref{typeI/IIsp} include
fractionally integrated class $|1- \e^{- \i x}|^{-2d_1}  |1- \e^{- \i y}|^{-2d_2}$ discussed in
\cite{bois2005}, \cite{lav2007}, \cite{guo2009}. Notice that $f_{\rm I}$ has a unique singularity at $(0,0)$ while
$f_{\rm II}$ is singular  on both coordinate axes and factorizes at low frequencies into a product of two functions depending on $x$ and $y$ alone.  See
Fig.1 below.
\vskip.1cm

 \medskip

\begin{figure}[ht!]
\begin{tabular}{cc}
\includegraphics[width=0.40\textwidth]{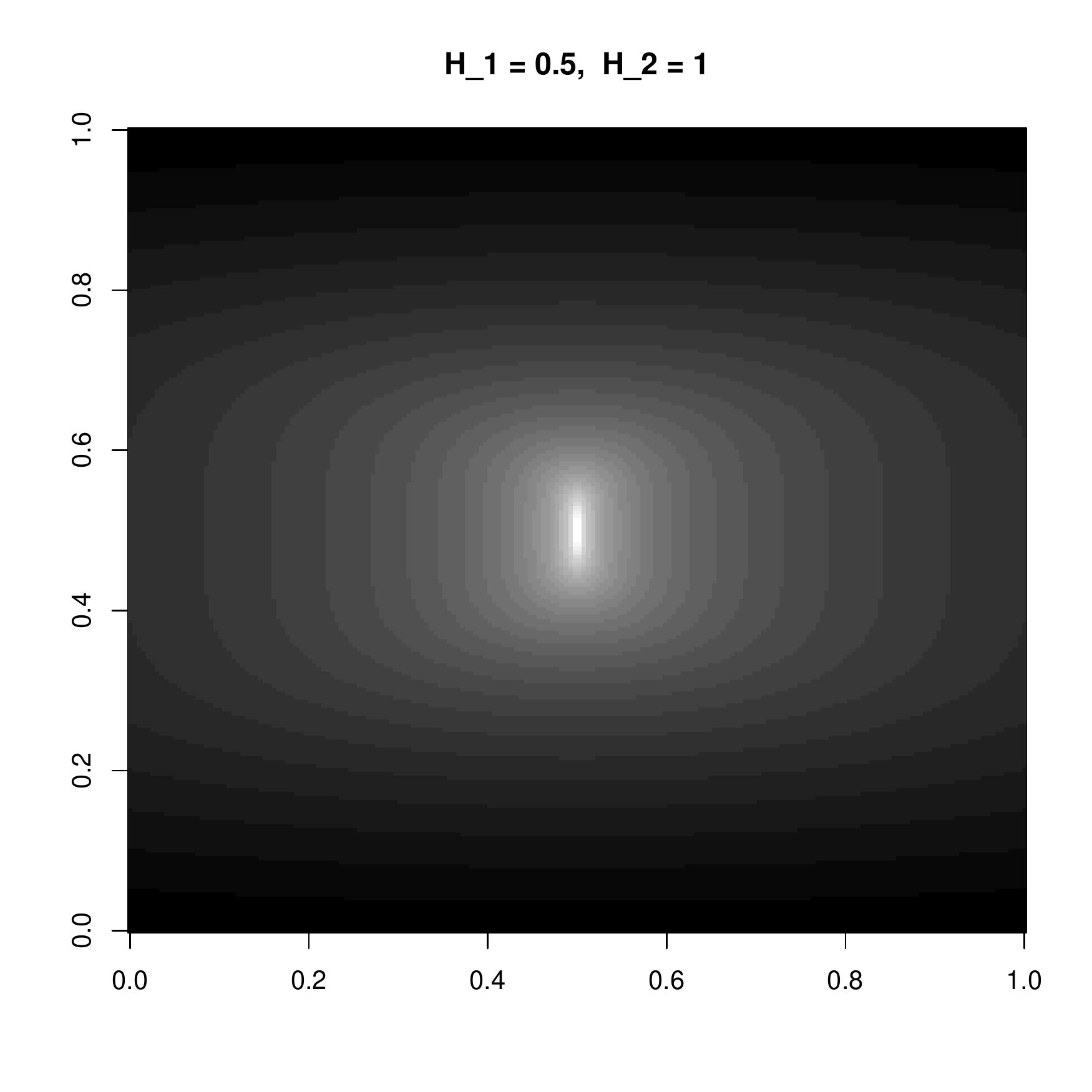}
&\includegraphics[width=0.40\textwidth]{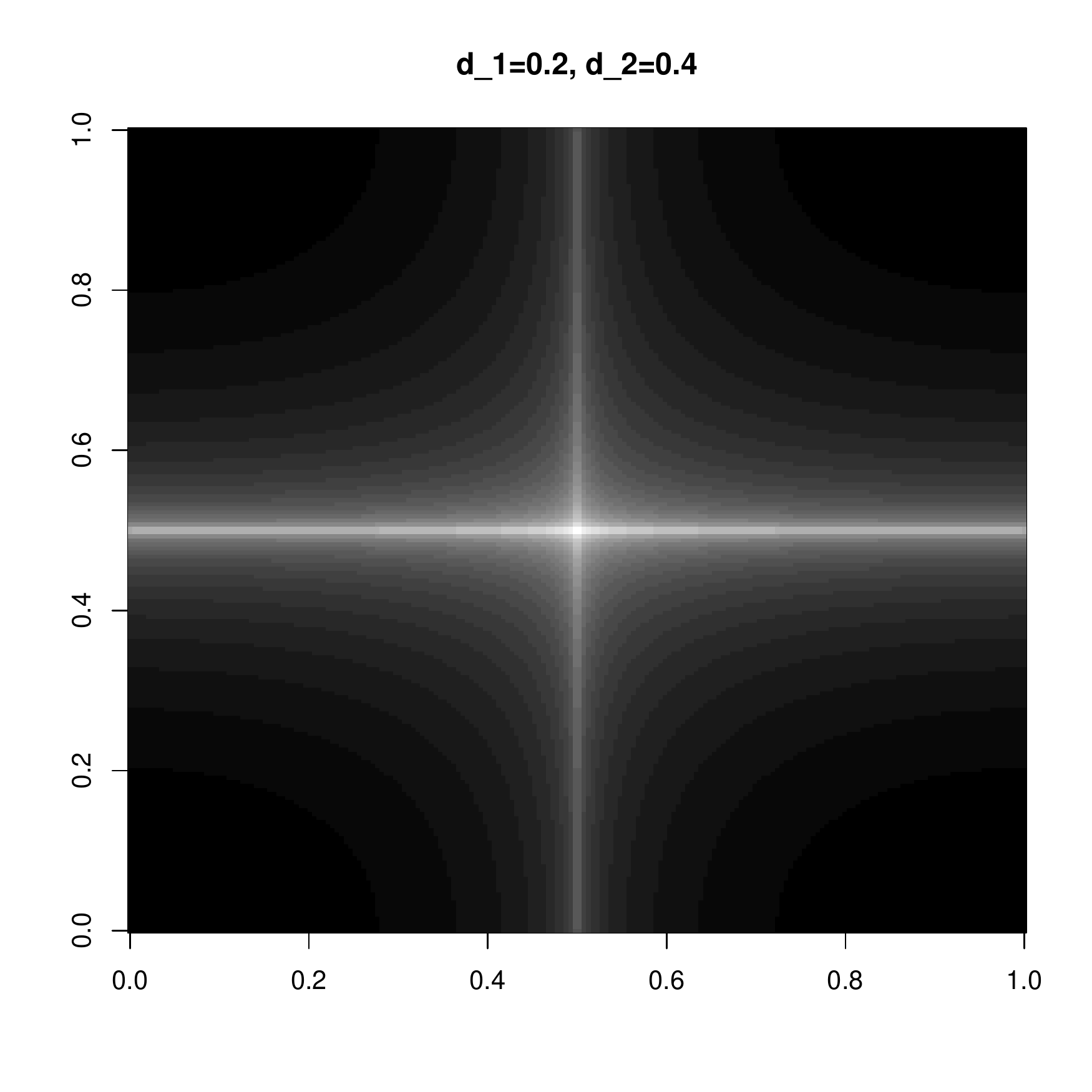}
 \\
{\small Type I spectral density $f_{\rm I}$ of \eqref{typeI/IIsp}, $H_1 = 0.5, H_2 = c=1$} &
{\small Type II spectral density  $f_{\rm II}$ of \eqref{typeI/IIsp}, $d_1 = 0.2, d_2 = 0.4$}
\\
\end{tabular}
\end{figure}

\smallskip

\centerline{Figure 1}

\bigskip


\smallskip

The main result of the present paper is Theorem \ref{thmGauss} which says that for Gaussian RFs with spectral density $f_{\rm I}$,
scaling transition occurs at $\gamma_0 = H_1/H_2$. It turns out that for such RFs  the unbalanced    
scaling limits
$V_+ $ and $V_-$  agree, up to a multiplicative constant, with a fractional Brownian sheet $B_{{\cal H}_1, {\cal H}_2}$
where at least one of  the two parameters
${\cal H}_1, {\cal H}_2 $ equals $1/2$ or 1. Recall that a fractional Brownian sheet $B_{{\cal H}_1, {\cal H}_2}$ with parameters
$0< {\cal H}_1, {\cal H}_2 \le 1$
is a Gaussian process on $\bar \R^2_+$ with
zero mean and covariance function
\begin{equation}\label{FBs}
\E B_{{\cal H}_1,{\cal H}_2}(x,y) B_{{\cal H}_1,{\cal H}_2}(x',y') = (1/4) (x^{2{\cal H}_1} + x'^{2{\cal H}_1} - |x-x'|^{2{\cal H}_1})
 (y^{2{\cal H}_2} + y'^{2{\cal H}_2} - |y-y'|^{2{\cal H}_2}),
\end{equation}
$(x,y), (x',y') \in \bar \R^2_+$. Particularly, for ${\cal H}_1 =1/2$,
$B_{{\cal H}_1, {\cal H}_2}(x,y)$ is a usual Brownian motion in $x$  having independent increments
in the horizontal direction and, for ${\cal H}_1 =1$, $ B_{{\cal H}_1, {\cal H}_2}(x,y)
= x B_{{\cal H}_2}(y) $ is random line in $x$ having shift-invariant (completely dependent) increments in the horizontal
direction  (see sec.2 for the definitions). The case when ${\cal H}_2$
equals 1/2 or 1 is analogous.
One may conclude that the unbalanced 
limits of the above Gaussian RFs
have a very special dependence structure (either independence or extreme (`deterministic') dependence along one of the coordinate axes).
By contrast, the well-balanced scaling limit $V_{\gamma_0}$ is not a fractional Brownian sheet and
has dependent but not shift-invariant increments
in arbitrary direction  on the plane. The dependence properties of rectangular increments
are made formal in sec.2 leading to the notion of Type I distributional LRD and isotropic/anisotropic LRD properties
for RFs on $\Z^2$. As shown in Proposition \ref{farimaLM},
stationary Gaussian RFs with spectral density $ f_{\rm II}$ in \eqref{typeI/IIsp} do not
exhibit scaling transition since in this case, all scaling limits $V_\gamma, \gamma >0 $ are equal to
$B_{{\cal H}_1, {\cal H}_2}$ with ${\cal H}_i = d_i + (1/2), i=1,2$ up to a multiplicative constant.

The above mentioned differences in the scaling behavior of Gaussian RFs with spectral densities
 $ f_{\rm I}$ and  $ f_{\rm II}$  \eqref{typeI/IIsp} are reflected in the scaling behavior of these spectral
densities. Indeed, the point $\gamma_0  = H_1/H_2$ at which  scaling  transition occurs  in the case of $f_{\rm I}$
can be characterized as a unique point $\gamma = \gamma_0 >0$ for which
a `non-degenerated' limit
\begin{equation}
\lim_{\lambda \to 0}\lambda^{H_1}f_{\rm I}(\lambda x, \lambda^\gamma y)  =   \big(|x|^2 + c|y|^{2 H_2/H_1}\big)^{-H_1/2}
\label{fIlim}
\end{equation}
exists, since for $\gamma \neq  \gamma_0$ the limit in \eqref{fIlim}
is either zero  or
`degenerated', in the sense that it does not depend on $y$. On the other hand, in the case of  $f_{\rm II}$,
a  `non-degenerated'  scaling limit $\lim_{\lambda \to 0}\lambda^{2d_1+ 2d_2 \gamma}f_{\rm II}(\lambda x, \lambda^\gamma y)\ =\
|x|^{-2d_1}|y|^{-2d_2} $ exists for any $\gamma >0$ and does not depend on $\gamma $.

It is of interest to extend
the results of this paper in several directions. Paper
\cite{lls2014} obtains scaling limits of LRD Gaussian RFs with singular spectral density
$f(x,y) =  g(x,y)   \big(|x- \mu y|^2 + c|y|^{2 H_2/H_1}\big)^{-H_1/2}$
having a general anisotropy axis $x - \mu y = 0,  \mu \in \R$ instead of $x =0$ in   $ f_{\rm I}$.
Further possibilities include investigation of scaling limits in  \eqref{Xsum01} for
{\it nonlinear} instantaneous functions $X(t,s) = G(Y(t,s))$
of stationary Gaussian RFs $Y = \{Y(t,s); (t,s) \in \Z^2\}$ with spectral density $f_{\rm I}$ in  \eqref{typeI/IIsp}
and
{\it non-Gaussian} moving-average RFs
\begin{equation}\label{Xlin}
X(t,s) = \sum_{(u,v) \in \Z^2} a(t-u, s-v) \vep(u,v), \qquad (t,s) \in \Z^2,
\end{equation}
where $\{ \vep(u,v); (u,v) \in \Z^2\}$ is an i.i.d. sequence with zero mean and finite variance, and $a(t,s)$ are deterministic
coefficients having the form
\begin{equation}
a(t,s) = \frac{g(t,s)}{ (|t|^2 + |s|^{2q_2/q_1})^{q_1/2}},  \qquad (t,s) \in \Z^2,
\end{equation}
where $g(t,s), (t,s) \in \Z^2 $ are bounded with $\lim_{|t|+|s| \to \infty} g(t,s) = 1$ and
$q_1, q_2 >0$ satisfy $(q_1 + q_2)/2 < q_1 q_2 < q_1 + q_2$. These  conditions
guarantee that $\sum_{(t,s) \in \Z^2} |a(t,s)|^2 < \infty,
\sum_{(t,s) \in \Z^2} |a(t,s)| = \infty $, hence \eqref{Xlin} is a well-defined LRD RF. We conjecture that RF $X$ in
 \eqref{Xlin} exhibits scaling transition at $\gamma_0 = q_1/q_2$ with
$V_+, V_-, V_{\gamma_0}$ similar as in Theorem  \ref{thmGauss} and $H_1, H_2 $ related to $q_1, q_2 $ by
$H_1 = (2/q_2)(q_1 + q_2 - q_1 q_2),  H_2 = (2/q_1)(q_1 + q_2 - q_1 q_2)$.
See  \cite{miko2002}, \cite{gaig2003} and  Remark \ref{remON} on a different type of scaling transition in telecommunication models.
A challenging task is  generalization of our limit results  to  sums $S_{n {\small \mbgamma}}({\mbf x})=
\sum_{{\mbf t} \in \prod_{i=1}^\nu [1, n^{\gamma_i} x_i]}
X({\mbf t}),\, {\mbf x} = (x_1, \dots, x_\nu) \in \R^\nu_+, \,
\mbgamma = (\gamma_1, \dots, \gamma_\nu) \in \R^\nu_+$ of
stationary Gaussian or linear RFs
on $\Z^\nu, \, \nu > 2$ with spectral
density $f({\mbf x}), {\mbf x}  \in [-\pi,\pi]^\nu $
similar to  $f_{\rm I}$ in  \eqref{typeI/IIsp} and
having $\nu$ parameters $H_1>0, \dots, H_\nu >0$. We note that,  instead of the single `balance condition'  $\gamma_2/\gamma_1 = \gamma_0 = H_1/H_2$ when $\nu=2$,
in higher dimensions $\nu >2$ there are $\nu (\nu-1)/2 >1$ `balance
conditions'   $\gamma_i/\gamma_j = H_j/H_i, \, i\ne j, 1\le i,j \le \nu $.
Depending on which of these `balance conditions'
are fulfilled or violated, we may expect
different scaling limits of  $S_{n {\small \mbgamma}}(\mbx)$.
We plan to explore the case $\nu=3$ in a forthcoming paper.


The notion of scaling transition for RFs is intrinsically related to the LRD property, which is often
identified with unboundedness of spectral density. It is clear that an i.i.d.  RF with zero mean and finite variance
does not exhibit  scaling transition
since its all scaling limits $V_\gamma, \gamma >0$ agree with Brownian sheet $B_{1/2,1/2}$. A similar fact remains true
for weakly dependent stationary RFs satisfying some mixing or other weak dependence  conditions.  On the other hand,
scaling limits of LRD RFs form a very rich class and are extensively studied. See, e.g.,  \cite{alb1994},
\cite{AnhLRM2012},
\cite{dobmaj1979}, \cite{douk2002},
\cite{douk2003p},  \cite{lav2007}, \cite{leo1999}, \cite{leool2012}, \cite{leo2011},
\cite{ps2014}, \cite{sur1982} and the references therein. Stationary Gaussian RFs form probably the most simple class of
LRD RFs, for which the asymptotic scaling theory is well-developed \cite{dob1979}. Nevertheless, we think that
our results  shed a new angle  on Gaussian RFs and spatial LRD.

\vskip.2cm

{\it Notation.} In what follows, $C$ denotes a generic constant
which may be different at different locations. We write $\limfdd,    \eqfdd, $ and $ \neqfdd $  for the weak convergence, equality and inequality
of finite-dimensional
distributions, respectively.
$\R^\nu_+ := \{ (x_1, \dots, x_\nu) \in \R^\nu: x_i > 0,
i=1, \dots, \nu \},  \bar \R^\nu_+ := \{ (x_1, \dots, x_\nu) \in \R^\nu: x_i \ge 0,
i=1, \dots, \nu \}, \,
\R_+ := \R^1_+, \bar \R_+ := \bar \R_+, \, \R^2_0 := \R^2 \setminus \{(0,0)\}$.
$\1(A)$ stands for the  indicator function of
a set $A$.
All equalities and inequalities between
random variables are assumed to hold almost surely.

\section{Distributional LRD properties of RFs on $\Z^2$ }

This sec. presents the definition of Type I distributional LDR property for RFs on $\Z^2$ and some related definitions
introduced in our paper \cite{ps2014}.
It is well-known that partial sum limits, or scaling limits, characterize dependence properties of
random processes indexed by $\Z$. Particularly, let $X = \{ X(t); t \in \Z\}$ be a stationary process such that its partial sums
tend to a process $V = \{V(x); x \ge 0\}$ in the sense that
$n^{-H} \sum_{t=1}^{[nx]} X(t) \limfdd  V(x)$ for some $H >0$. If the limit process $V$ has dependent increments,
the process $X$ is  said to have {\it distributional long memory}, or {\it distributional LRD} property. The last property
(originating to Cox \cite{cox1984})
was introduced in Dehling and Phillips \cite{deh2002} and later used and verified for several classes of LRD processes
\cite{lps2005}, \cite{gir2009}, \cite{pps2014}, \cite{ps2010}.

Let $\ell = \{ (x,y) \in \R^2: a x + b y = c \} $ be a line in $\R^2. $
A line
$\ell' = \{ (x,y) \in \R^2: a' x + b' y = c' \} $ is said
{\it perpendicular to $\ell $} (denoted  $\ell ' \bot \ell $) if $a a'+ b b' = 0$. Write $(u,v) \prec (x,y)$ (respectively,
$(u,v) \preceq (x,y)$), $(u,v), (x,y)  \in \R^2$
if $u < x$ and $v < y$ (respectively, $u \le x$ and  $v \le y$) hold.
A {\it rectangle}
is a set  $K_ {(u,v); (x,y)} := \{ (s,t) \in \R^2:   (u,v) \prec (s,t) \preceq (x,y) \};  K_{x,y} := K_{(0,0); (x,y)}.    $
Denote $K_ {(u,v); (x,y)} + (z,w) := K_ {(u+z,v+w); (x+z,y+w)}$ the rectangle $K_ {(u,v); (x,y)}$ shifted by $(z,w) \in \R^2$.
We say that two
rectangles $K = K_{(u,v); (x,y)} $ and $K'= K_{(u',v'); (x',y')} $
are {\it separated by line $\ell'  $} if they lie on different sides
of $\ell' $, in which case $K$ and $K'$ are necessarily disjoint: $K \cap K' = \emptyset $.
See Fig.~2.

\begin{center}
\begin{figure}[h]
\begin{center}
\includegraphics[width=8 cm,height=5.5cm]{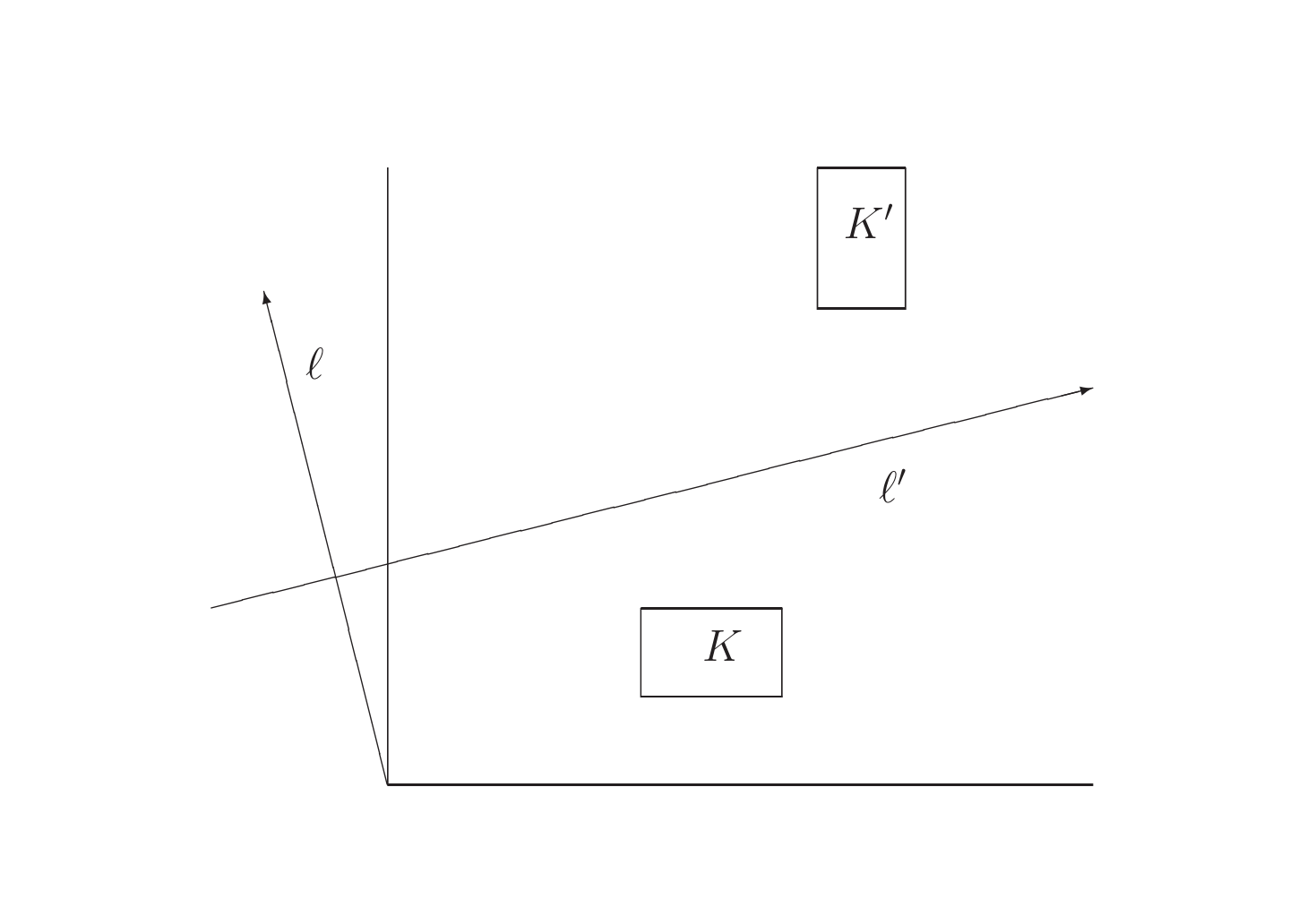}
\end{center}
\end{figure}
\end{center}

\smallskip

\centerline{Figure 2}

\bigskip

Let $V =  \{V(x,y); (x,y) \in \bar \R^2_+ \}$ be a
RF and $K = K_{(u,v); (x,y)} \subset \R^2_+ $ be a rectangle.
By {\it increment of $V$ on rectangle $K$}  we mean the difference
$$
V(K) :=  V(x,y) - V(u,y) - V(x,v) + V(u,v).
$$
We say that $V$ {\it  has stationary rectangular increments} if
\begin{equation} \label{Vstatinc}
\{ V(K_ {(u,v); (x,y)}); (u,v) \preceq (x,y)\}  \eqfdd \{V(K_ {(0,0); (x-u,y-v)});
(u,v) \preceq (x,y) \}, \qquad
\text{for any} \quad (u,v)  \in \R^2_+.
\end{equation}

\begin{defn} \label{iinc}
Let $V= \{V(x,y); (x,y) \in \bar \R^2_+ \}$ be a RF with stationary rectangular increments,
$V(x,0) = V(0,y) \equiv 0, \, x,y \ge 0,$ and $\ell \subset \R^2$ be a given line , $(0,0) \in \ell $.
We say that $V$ has

\smallskip

\noi (i) {\rm  independent rectangular
increments
in direction $\ell $ } if for any orthogonal line $\ell ' \bot  \ell$ and any two rectangles $K, K' \subset \R^2_+$ separated by $\ell '$,
increments
$V(K)$ and $ V(K')$
are independent;

\smallskip

\noi (ii) {\rm  invariant rectangular
increments
in direction $\ell $ } if $V(K) = V(K')$
for any two rectangles  $K, K' \subset \R^2_+$ such that
$K' = (x,y) + K$ for some $(x,y) \in \ell $;

\smallskip

\noi (iii) {\rm dependent rectangular
increments
in direction $\ell $ } if neither (i) nor (ii) holds;

\smallskip

\noi (iv) {\rm dependent rectangular
increments} if   $V$ has dependent rectangular increments
in arbitrary direction;

\smallskip

\noi (v) {\rm  independent rectangular
increments} if  $V$     has independent rectangular increments
in arbitrary direction.

\end{defn}

\begin{ex} \label{exFBS}
{\rm {\it Fractional Brownian sheet} $B_{{\cH}_1,{\cH}_2}$ with parameters  $0 < {\cH}_1, {\cH}_2 \le 1$ is a
Gaussian process on $\bar \R^2_+$ with zero mean and covariance in \eqref{FBs}.
It follows (see \cite{aya2002}, Cor.3)
that for any rectangles
$K = K_{(u,v); (x,y)}, K'  = K_{(u',v'); (x',y')}$
\begin{eqnarray}
&&\E B_{{\cH}_1,{\cH}_2}(K) B_{{\cH}_1,{\cH}_2}(K') \nn \\
&&= \ \E(B_{{\cH}_1}(x)-B_{{\cH}_1}(u))(B_{{\cH}_1}(x')-B_{{\cH}_1}(u')) \E (B_{{\cH}_2}(y)-B_{{\cH}_2}(v))(B_{{\cH}_2}(y')-B_{{\cH}_2}(v')),  \label{FBsheet}
\end{eqnarray}
where $\{B_{{\cH}_i}(x); x\in \bar\R_+\}$ is a fractional Brownian motion on $\bar \R_+ = [0,\infty)$ with
$\E B_{{\cH}_i}(x) B_{{\cH}_i}(x') = (1/2)(x^{2{\cH}_i} + x'^{2{\cH}_i} - |x-x'|^{2{\cH}_i}, \, i=1,2.$ (For ${\cH} = 1$, the process
 $\{B_{{\cH}}(x); x\in \bar \R_+\}$ is a random line.)
In particular, $B_{{\cH}_1,{\cH}_2}$ has stationary rectangular increments, see (\cite{aya2002}, Prop.2). It follows from (\ref{FBsheet})
(see \cite{ps2014} for details) that
fractional Brownian sheet $B_{{\cH}_1,{\cH}_2}$ has:
\begin{itemize}
\item dependent rectangular increments if ${\cH}_i \not\in \{1/2, 1\}, i=1,2$;

\item independent rectangular increments in the horizontal (vertical) direction
if  ${\cH}_1 = 1/2 $ (${\cH}_2 = 1/2 $);

\item invariant rectangular increments in the horizontal (vertical) direction
if  ${\cH}_1 = 1 $ (${\cH}_2 = 1$);

\item independent rectangular increments if ${\cH}_1 = {\cH}_2 = 1/2 $.

\end{itemize}}
\end{ex}

\begin{defn} \label{Adlm} Let $X = \{X(t,s); (t,s) \in \Z^2 \}$ be
a stationary RF. Assume that for any $\gamma >0$ there  exist normalization
$A_n(\gamma) \to \infty $ and a RF \, $V_\gamma = \{V_\gamma(x,y); (x,y) \in \bar \R^2_+\}$  such that \eqref{Xsum01} holds.
We say that $X$ {\rm has Type I distributional LRD} (or {\rm  $X$ is a Type I RF})
if there exists a unique $\gamma_0 \in (0, \infty) $ such that:

\begin{itemize}

\item RF $V_{\gamma_0}$ has dependent rectangular increments; and

\item RFs  $V_\gamma, \gamma \neq \gamma_0$  do not have dependent rectangular increments; in other words,
for each $\gamma  \neq \gamma_0, \gamma >0$ there exists a line $\ell(\gamma)$ such that RF
$V_\gamma $ has either independent, or invariant rectangular increments in direction $\ell(\gamma)$.

\end{itemize}
Moreover, a Type I RF \, $X$ is said to have {\rm isotropic distributional LRD} if $\gamma_0 = 1$ and
{\rm anisotropic distributional LRD} if $\gamma_0 \neq 1$.

\end{defn}

\begin{rem}\label{Vanis} {\rm   Let $X = \{X(t,s); (t,s) \in \Z^2\}$ be a stationary RF satisfying \eqref{Xsum01} for some $\gamma >0$.
Then the limit RF $ V_\gamma $ has stationary rectangular increments in the sense of
\eqref{Vstatinc}. Moreover, if $A_n(\gamma) = n^{H}$ for some $H >0$, then   $ V_\gamma $
satisfies the following self-similarity property (see \cite{ps2014}):
\begin{equation}\label{OS1}
\{ \lambda V_\gamma(x,y); (x,y) \in \R^2 \}   \ \eqfdd \     \{ V_\gamma(\lambda^{1/H} x, \lambda^{\gamma/H} y); \, (x,y) \in \R^2 \},
 \quad \forall \, \lambda > 0.
\end{equation}
We note that \eqref{OS1}  is a particular case of {\it operator scaling} property for RFs
introduced in \cite{bier2007}. }

\end{rem}

\section{Main results}

This section obtains the presence/absence of scaling  transition for  stationary Gaussian
RFs with spectral densities in \eqref{typeI/IIsp}, and a characterization of the dependence properties
of their scaling limits. For a given a stationary RF $X = \{ X(t,s); (t,s) \in \Z^2\}$
denote
\begin{equation}\label{Sgamma}
S_{n\gamma}(x,y):= \sum_{(t,s) \in K_{[nx, n^\gamma y] }} X(t,s), \qquad (x,y) \in \R^2_+.
\end{equation}

First, consider a Type  I spectral  density $f = f_{\rm I}$ of \eqref{typeI/IIsp}:
\begin{equation}\label{typeIsp}
f(x,y) = \frac{g(x,y)}{\big(|x|^2 + c|y|^{2 H_2/H_1}\big)^{H_1/2}}, \quad (x,y) \in \Pi^2 = [-\pi,\pi]^2,
\end{equation}
where   $0< H_1 \le H_2 < \infty, H_1H_2 <H_1 + H_2,  c>0$ and $g$ is  bounded and continuous at the origin with $g(0,0)=1$.  We have
\begin{equation}  \label{h}
h(x,y)  := \lim_{\lambda \to 0} \lambda f(\lambda^{1/H_1}x, \lambda^{1/H_2}y) = \frac{1}{\big(|x|^2 + c|y|^{2 H_2/H_1}\big)^{H_1/2}}, \quad (x,y) \in \R^2_0.
\end{equation}
Note that $h$ is continuous on $\R^2_0$ and
satisfies the scaling property: for any $\lambda >0$
\begin{eqnarray}\label{hlimit}
\lambda h(\lambda^{1/H_1}x, \lambda^{1/H_2}y)&=&h(x,y), \qquad \forall (x,y) \in \R^2_0.
\end{eqnarray}
With $h$ in \eqref{h}
we associate a family of
Gaussian RFs indexed by $\gamma >0$, as follows. For $\gamma = \gamma_0 := H_1/H_2$, set
\begin{equation}  \label{Vlimit00}
V_{\gamma_0}(x,y) := \int_{\R^2}\frac{ (1-\e^{{\i} u x})(1 - \e^{{\i} v y})}{\i^2 u v}\, \sqrt{h(u,v)}   W(\d u, \d v),
\quad (x,y) \in \bar \R^2_+,
\end{equation}
where $\{W(\d x, \d y); \ (x,y) \in \R^2
\}$ is a standard complex-valued Gaussian noise, $\overline{W(\d x, \d y)} =
W(-\d x, - \d y)$, with zero mean and
variance $\E |W(\d x, \d y)|^2 = \d x \d y$.
Define
\begin{eqnarray}\label{Vlimit11}
V_+(x,,y)
&:=& \begin{cases}
\int_{\R^2}\frac{ (1-\e^{{\i} u x})(1 - \e^{{\i} v y})}{\i^2 u v}\, |u|^{-H_1/2}   W(\d u, \d v),
&H_1 <1, \\
x \rho_1 c^{(1-H_1)/4} \int_{\R}\frac{1-\e^{{\i} v y}}{\i v}\, |v|^{-(H_1H_2 -H_2)/2H_1} W_1(\d v),
&H_1>1,
\end{cases}
\end{eqnarray}
and
\begin{eqnarray}  \label{Vlimit12}
V_-(x,y)
&:=& \begin{cases}
\int_{\R^2}\frac{ (1-\e^{{\i} u x})(1 - \e^{{\i} v y})}{\i^2 u v}\, |v|^{-H_2/2}   W(\d u, \d v),
&H_2 < 1, \\
y \rho_2 c^{-H_1/4H_2}
\int_{\R}\frac{1-\e^{{\i} u x}}{\i u}\, |u|^{-(H_1H_2 -H_1)/2H_2} W_1(\d u),
&H_2>1.
\end{cases}
\end{eqnarray}
where  $\{W(\d x, \d y)
\}$ is as in \eqref{Vlimit00} and $\{W_1(\d x); \ x\in \R \}$ is
a standard complex-valued Gaussian noise on $\R$, $\overline{W_1(\d x)} =
W_1(-\d x)$, with zero mean and
variance $\E |W_1(\d x)|^2 = \d x, \, \rho^2_1 := B(1/2,(H_1-1)/2),  \rho^2_2 := (H_1/2H_2) B(H_1/H_2, (H_1H_2 - H_1)/2H_2).  $
Here and below, $B(\cdot, \cdot) $ is the beta function.
We also define
\begin{eqnarray}
H_+(\gamma)&:=&
\begin{cases}(1 + \gamma + H_1)/2, &H_1  < 1, \\
(\gamma H_1 + \gamma H_1 H_2 - \gamma H_2 + 2 H_1)/2H_1, &H_1 > 1,
\end{cases} \label{H+} \\
H_-(\gamma)
&:=&
\begin{cases}
(1 + \gamma + \gamma H_2)/2, &\quad H_2 < 1, \\
(H_2 + H_1 H_2 - H_1 + 2\gamma H_2)/2H_2, &\quad H_2 > 1,
\end{cases} \label{H-}
\end{eqnarray}
where $\gamma >0$, and
\begin{eqnarray}\label{H0}
H(\gamma_0)&:=&(H_1 + H_2 + H_1 H_2)/2H_2,  \hskip1cm H_1\ne 1, \ H_2 \ne  1.
\end{eqnarray}
Note $H_+(\gamma_0) = H_-(\gamma_0) = H(\gamma_0)$.

\begin{proposition}\label{Vgamma}
(i) The RFs $V_{\gamma_0} = \{V_{\gamma_0}(x,y); (x,y) \in \bar \R^2_+\}, $
$V_+  = \{V_{+}(x,y); (x,y) \in \bar \R^2_+\}$ and $V_-  = \{V_{-}(x,y); (x,y) \in \bar \R^2_+\}$
in  \eqref{Vlimit00}-\eqref{Vlimit12} are well-defined for any $0<H_1 \le H_2 <  \infty, \, H_1H_2 < H_1 + H_2, \, c>0$ with
exception of $H_1=1$ in \eqref{Vlimit11} and $H_2 =1 $ in \eqref{Vlimit12}.
These  RFs
have zero mean, finite variance
and stationary rectangular increments in the sense of \eqref{Vstatinc}. Furthermore,
\begin{eqnarray}\label{VB1}
V_+
&\eqfdd&\kappa^+_{H_1,H_2} \begin{cases}
B_{(H_1+1)/2, 1/2}\,
&H_1 <1, \\
B_{1, (1+ H_2- H_2/H_1)/2},
&H_1>1,
\end{cases} \\
V_-
&\eqfdd&\kappa^-_{H_1,H_2} \begin{cases}
B_{1/2, (H_2+1)/2}\,
&H_2 <1, \\
B_{(1+ H_1- H_1/H_2)/2, 1},
&H_2>1,
\end{cases} \label{VB2}
\end{eqnarray}
where $B_{H_1,H_2}$ is a fractional Brownian sheet (see Example  \ref{exFBS}), and $\kappa^\pm_{H_1, H_2} >0 $ are
some constants.

\medskip

\noi (ii) $V_{\gamma_0}, V_+, V_-$ are operator scaling RFs: for any $\lambda >0$,
\begin{eqnarray}
\{V_{\gamma_0}(\lambda x, \lambda^{\gamma_0} y); (x,y) \in \bar  \R^2_+\}
&\eqfdd& \{\lambda^{H(\gamma_0)} V_{\gamma_0}(x,y); (x,y) \in \bar \R^2_+\}, \label{V0}\\
\{V_{+}(\lambda x, \lambda^\gamma y); (x,y) \in \bar  \R^2_+\}&\eqfdd& \{\lambda^{H_+(\gamma)} V_{+}(x,y); (x,y) \in \bar \R^2_+\},
\label{V+}\\
\{V_{-}(\lambda x, \lambda^\gamma y); (x,y) \in \bar  \R^2_+\}&\eqfdd& \{\lambda^{H_-(\gamma)} V_{-}(x,y); (x,y) \in \bar \R^2_+\},\label{V-}
\end{eqnarray}
where  $H_+(\gamma), H_-(\gamma), H(\gamma_0) $ are defined in \eqref{H+}, \eqref{H-}, \eqref{H0}, respectively.
In \eqref{V+} and \eqref{V-},  $\gamma >0$ is arbitrary.

\smallskip

\noi (iii) $V_{\gamma_0}$ has dependent rectangular increments, while $V_+$ and $V_-$ have either independent, or invariant
rectangular increments along one of the coordinate axes.

\end{proposition}

\noi {\it Proof.} (i) Let us show that $V_{\gamma}(x,y)$ is well-defined as  stochastic integral
w.r.t. Gaussian white noise.
It suffices to consider the case $x=y=1$ only since the general  case is analogous.
Let $\gamma = \gamma_0$. Then
\begin{eqnarray*}  \label{Vvar}
\E V^2_{\gamma_0}(1,1)&=&\int_{\R^2} \big|1-\e^{{\i} u}\big|^2 \, \big|1 - \e^{{\i}v}\big|^2 \,
\frac{h(u,v)}{|u v|^2}\, \d u \d v   \ \le \ C\int_0^\infty  \int_0^\infty
\frac{\d u \d v}{(1+u^2)(1+ v^2)(u^2 + v^{2 H_2/H_1})^{H_1/2}} \\
&=&C\Big(\int_0^\infty \d u \int_0^1 \d v \dots + \int_0^\infty \d u \int_1^\infty \d v \dots \Big) =: C(J_1 + J_2).
\end{eqnarray*}
By change of variable $v = u^{H_1/H_2} z $,  $
J_1 \le C \int_0^\infty \frac{\d u}{1+ u^2} \int_0^{1/u^{H_1/H_2}}
\frac{ u^{H_1/H_2} \d z}{u^{H_1} (1 + z^{2 H_2/H_1})^{H_1/2}} =: J'_1.    $
Let $H_2 > 1$. Then $J_1' \le C \int_0^\infty u^{(H_1/H_2) - H_1} (1+u^2)^{-1} \d u  < \infty $ since
$H_2(H_1 -1) < H_1. $  Next,
let $H_2  < 1$. Then  $J'_1 \le C \int_0^\infty (1+u^2)^{-1} \d u  < \infty $ and
$J_1 < \infty$. The case $H_1 =1 $ follows similarly.
The convergence $J_2 < \infty$ is obvious.

Let us show that \eqref{Vlimit11} is well-defined. Let $H_1 < 1$, then
$\E V^2_{\gamma}(1,1) \le C\int_{\R^2_+}  \d u \d v (1+u^2)^{-1}(1+ v^2)^{-1} v^{-H_1} < \infty $. Next, let $H_1 > 1$, then
$\E V^2_{\gamma}(1,1) \le C\int_0^\infty  \d v (1+ v^2)^{-1} v^{-(H_1H_2 - H_2)/H_1} < \infty $ since $H_2 (H_1 -1) <  H_1$. The convergence
of the stochastic integral in \eqref{Vlimit12} follows in a similar way.

Let  $K_{(x,y); (x',y')}, (x,y) \preceq (x',y') $ be a rectangle in $\R^2_+$. Then from
\eqref{Vlimit00} we immediately obtain $ V_{\gamma_0}(K_ {(x,y); (x',y')}) $   $ =   \tilde V_{\gamma_0} (K_ {(0,0); (x'-x,y'-y)})$, where
$\tilde V_{\gamma_0}(x,y) $ is defined as in \eqref{Vlimit00} with $W(\d u, \d v)$ replaced by
$\tilde W(\d u, \d v) :=  \e^{\i (ux + vy)} W(\d u, \d v). $ Clearly $\{\tilde W  (\d u, \d v); (u,v) \in \R^2 \}
\eqfdd \{W  (\d u, \d v); (u,v) \in \R^2 \}$ for any $(x,y) \in \R^2_+$. Hence, $V_{\gamma_0} $ has stationary
rectangular increments. The same fact for $V_\pm$ follows analogously.

Relation \eqref{VB1} follows from Gaussianity and
\begin{eqnarray}
&&\E [V_+ (x,y) V_+(x',y')] \nn \\
&&\ = \ (\kappa^+_{H_1,H_2})^2
\begin{cases}
\E [B_{(H_1+1)/2}(x)B_{(H_1+1)/2}(x')]  \E [B_{1/2}(y)B_{1/2}(y')], &H_1 < 1, \\
 \E [B_{1}(x)B_{1}(x')] \E[  B_{(1+ H_2- H_2/H_1)/2}(y) B_{(1+ H_2- H_2/H_1)/2}(y')],
&H_1 >  1
\end{cases}\label{Vfbm}
\end{eqnarray}
for any $x, x', y, y' \ge 0$. Let $H_1 < 1$, then the l.h.s.  of \eqref{Vfbm} factorizes as the product of two integrals
$\int_{\R} (1- \e^{\i u x})(1- \e^{-\i ux'}) |u|^{-2- H_1} \d u \,  \int_{\R} (1- \e^{\i v y})(1- \e^{-\i vy'}) |v|^{-2} \d v $,
equal to the covariances on the r.h.s., see e.g. \cite{taq2003}, Prop.9.2. The case  $H_1 > 1$ in   \eqref{Vfbm}
and \eqref{VB2} are analogous.

\smallskip

\noi (ii)   The operator scaling property follows from scaling properties of the integrands, see
\eqref{hlimit}, and the white noise, viz.,
$\{W(\d u/\lambda, \d v/\lambda^\gamma); (u,v) \in \R^2\}
\eqfdd \{\lambda^{-(1+\gamma)/2} W(\d u, \d v); (u,v) \in \R^2 \}, \ \{W_1(\d u/\lambda); u \in \R\}
\eqfdd \{\lambda^{-1/2} W_1(\d u); u \in \R \}. $

\smallskip

\noi (iii) The fact that $V_+, V_- $ do not have dependent rectangular increments follows from  \eqref{VB1}, \eqref{VB2}
and the properties of fractional Brownian sheet stated in Example \ref{exFBS}.
The proof that $V_{\gamma_0}$ in \eqref{Vlimit00}  has dependent increments (i.e., neither independent nor invariant rectangular increments in any direction)
is part of a more general  statement in Lemma \ref{lemspec} below.
Proposition \ref{Vgamma}  is proved.  \hfill $\Box$


\begin{theorem} \label{thmGauss} Let $X$
be a stationary zero-mean Gaussian RF on $\Z^2$ with zero mean and spectral
density  $f$  in  \eqref{typeIsp}, where
$c>0, \, 0 < H_1 \le H_2<\infty, H_1, H_2 \neq 1 $ and $H_1 H_2 < H_1 + H_2  $.
Then for any $\gamma >0$ the limit of partial sums
\begin{equation} \label{TypeIlim}
n^{-H(\gamma)} S_{n\gamma} (x,y)
\limfdd \ V_\gamma (x,y), \quad (x,y) \in  \R^2_+, \quad n \to \infty
\end{equation}
exists with
\begin{equation}
V_\gamma\  :=\  \begin{cases}V_+, &\gamma > \gamma_0, \\
V_-, &\gamma < \gamma_0, \\
V_{\gamma_0}, &\gamma= \gamma_0,
\end{cases}, \qquad
H(\gamma)\ :=\  \begin{cases}H_+(\gamma), &\gamma > \gamma_0, \\
H_-(\gamma), &\gamma < \gamma_0, \\
H(\gamma_0), &\gamma= \gamma_0,
\end{cases}
\end{equation}
and $V_\pm, V_{\gamma_0}, H_\pm(\gamma), H(\gamma_0)$  given in
\eqref{Vlimit00}-\eqref{Vlimit12} and
 \eqref{H+}-\eqref{H0}, respectively.
As a consequence, the RF $X$ exhibits  scaling transition at $\gamma_0 = H_1/H_2$.
Moreover, \ $ X$ has Type I isotropic distributional LRD if $H_1 = H_2$ and
Type I anisotropic distributional LRD if $H_1 \neq H_2$.
\end{theorem}

\begin{rem} {\rm The existence of the limit in \eqref{TypeIlim} in the cases $\gamma > \gamma_0, H_1 = 1$ and
$\gamma < \gamma_0, H_2  = 1$ is
an open question. Note that $V_+ $ in \eqref{Vlimit11} is undefined for $H_1 =1 $ and, similarly, $V_- $ in  \eqref{Vlimit12}
is undefined for $H_2 =1$. Nevertheless, relations \eqref{VB1}-\eqref{VB2} suggest
that the limit   \eqref{TypeIlim} might  exist also in the above cases and be given by $V_+ = \kappa_+ B_{1, 1/2} $
and $V_- = \kappa_- B_{1/2, 1} $ for some constants $\kappa_\pm \ne 0$.
}
\end{rem}

\noi {\it Proof of Theorem \ref{thmGauss}.} Let us prove the convergence in (\ref{TypeIlim}).
Recall the definition of $S_{n \gamma}(x,y) $ in \eqref{Sgamma}.
By Gaussianity,
this  follows from
\begin{equation} \label{Hlimvar}
R_{n \gamma}(x,y; x',y') := n^{-2H(\gamma)} \E [S_{n\gamma}(x,y) S_{n \gamma}(x',y')]  \ \to \ \E [V_\gamma(x,y) V_\gamma (x',y')], \qquad n \to \infty.
\end{equation}
We have, with $m := n^{\gamma},
$ and $D_n(u) := \sum_{t=1}^n  \e^{\i tu}  = (\e^{\i u}- \e^{\i (n+1)u})/(1- \e^{\i u}), \, |u| < \pi $,
\begin{eqnarray*}
R_{n \gamma}(x,y; x',y')
&=&n^{-2H(\gamma)} \int_{\Pi^2} D_{[nx]}(u) \overline{D_{[n x']}(u)} D_{[m y ]}(v)  \overline{D_{[m y']}(v)}  f(u,v) \d u \d v.
\end{eqnarray*}

Consider \eqref{Hlimvar} for   $\gamma = \gamma_0$. By change of variables, $R^\gamma_{n}(x,y; x',y')= \int_{\R^2} G_n(u,v) \d u \d v$,
where
\begin{equation}  \label{Gnuv}
G_n(u,v)\ :=\
\frac{1}{n^2 m^2 } D_{[nx]}(\mbox{$\frac{u}{n}$}) \overline{D_{[n x']}(\mbox{$\frac{u}{n}$})}
D_{[m y ]}(\mbox{$\frac{v}{m}$})  \overline{D_{[m y']}(\mbox{$\frac{v}{m}$})}
n^{-H_1}
f(\frac{u}{n}, \frac{v}{m}) \1(|u| \le \pi n,
|v| \le \pi m)
\end{equation}
and we used the fact that $n m n^{H_1} = n^{2H(\gamma_0)}$.
From \eqref{h}  it follows that
$n^{-H_1} f(\frac{u}{n}, \frac{v}{m}) \to h(u,v) $ a.e. in $\R^2$
and hence
$$
G_n(u,v)\ \to \  G(u,v) \ := \  \big(\frac{1- \e^{\i x u}}{\i u}\big)
\big(\frac{1- \e^{-\i x' u}}{-\i u}\big) \big( \frac{1- \e^{\i y v}}{\i v}\big) \big( \frac{1- \e^{-\i y' v}}{-\i v}\big)  h(u,v) \quad
\text{a.e. in} \quad  \R^2,
$$
where $\int_{\R^2} G(u,v) \d u \d v = \E [V_{\gamma_0}(x,y) V_{\gamma_0} (x',y')]
$. Moreover, $|n^{-1}D_{[nx]}(\mbox{$\frac{u}{n}$})| \le C x (1+ |([nx]/n) u|) \le C/(1+ |u|), |u|< n \pi $ for any fixed $x \in \R, x \ne 0$.
Together with $f(u,v) \le C h(u,v)$, see  \eqref{typeIsp}, \eqref{h},
this implies for any fixed $x,  x' , y, y' >0$ that
$$
|G_n(u,v)| \ \le \ C(1+ u^2)^{-1} (1+ v^2)^{-1} h(u,v)\  =: \ \bar G(u,v),
$$
where $\int_{\R^2} \bar G(u,v) \d v \d v < \infty $ (see above). Therefore, (\ref{Hlimvar}) for $\gamma =\gamma_0$ follows by the dominated convergence
theorem.

\smallskip

Consider (\ref{Hlimvar}) for $\gamma > \gamma_0, 0< H_1 < 1$. We have again
$R_{n \gamma}(x,y; x',y')= \int_{\R^2} G_n(u,v) \d u \d v$ with $G_n$ given in \eqref{Gnuv} and
$n^{-H_1} f(\frac{u}{n}, \frac{v}{m}) \to h(u,0) = |u|^{-H_1} $ a.e. in $\R^2$
and $n^{-H_1} f(\frac{u}{n}, \frac{v}{m}) \le C h(u,0), u \in \R$.
Since $\bar G(u,v) := C(1+ u^2)^{-1} (1+ v^2)^{-1} h(u,0)$ is  integrable on $\R^2$ for $0< H_1 < 1$, this proves
(\ref{Hlimvar}).

\smallskip

Consider (\ref{Hlimvar}) for $\gamma > \gamma_0, H_1 > 1$. Then
\begin{eqnarray*}
R_{n \gamma}(x,y; x',y')&=&\int_{\R^2} \tilde G_m(u,v) \d u \d v, \qquad \text{with} \\
\tilde G_m(u,v)&:=&L_{m1}(u) L_{m2}(v) L_{m3}(u,v) L_{4}(u,v)\1(|u| \le \pi m^{1/\gamma_0},
|v| \le \pi m),
\end{eqnarray*}
where $m = n^\gamma, \, m^{(1/\gamma)- (1/\gamma_0)} \to 0 $ and
\begin{eqnarray}\label{LLconv}
L_{m1}(u)&:=&\frac{n^{2((\gamma/\gamma_0)-1)}(1- \e^{\i u ( [nx]/n^{\gamma/\gamma_0})})(1- \e^{-\i u ( [nx']/n^{\gamma/\gamma_0})}) }
{|n^{\gamma/\gamma_0}(1- \e^{\i u/n^{\gamma/\gamma_0}})|^2} \ \to \   x x',  \\
L_{m2}(v)&:=&\frac{(1- \e^{\i v ([my]/m)})(1- \e^{-\i v ([my']/m)}) }
{|m (1- \e^{\i v/m})|^2} \ \to \ \frac{(1- \e^{\i v y})(1- \e^{- \i vy'})  }{v^2}, \nonumber \\
L_{m3}(u,v)&:=&g(\frac{u}{m^{1/\gamma_0}}, \frac{v}{m}) \ \to \  1, \nonumber \\
L_{4}(u,v)&:=&\frac{1}{(u^2 + cv^{2/\gamma_0})^{H_1/2}},   \nonumber
\end{eqnarray}
as $m \to  \infty$.
Note that for fixed $ x,x',y,y' >0$,
all three convergences in \eqref{LLconv} are uniform in $(u,v) \in \R^2 $  on each compact set in $\R^2$, the limit functions being
bounded and continuous in $\R^2$, moreover, $L_{mi}, i= 1,2,3 $ are bounded on  $|u| \le \pi m^{1/\gamma_0},
|v| \le \pi m $. Therefore,  by the dominated convergence theorem, as $n \to \infty $,
\begin{eqnarray*}
R_{n \gamma}(x,y; x',y')
&\to&xx'\int_{\R} \frac{(1- \e^{\i  v y})(1 - \e^{-\i v y'}) }{v^2} \d v \int_\R  \frac{\d u}{(u^2 + cv^{2/\gamma_0})^{H_1/2}} \\
&=&\rho_1^2 c^{(1-H_1)/2} x x' \int_{\R} \frac{(1- \e^{\i  v y})(1- \e^{-\i vy'})}{v^2} \frac{\d u}{|v|^{(H_1H_2 - H_2)/H_1}},
\end{eqnarray*}
where $\rho_1^2 = \int_\R  \frac{\d u}{(u^2 + 1)^{H_1/2}} = B(1/2, (H_1-1)/2))$. The last limit agrees with the covariance
on the r.h.s. of \eqref{Hlimvar}, see the definition of $V_+ = V_\gamma $ in \eqref{Vlimit11}, proving
(\ref{Hlimvar}) for $\gamma > \gamma_0$. The proof of (\ref{Hlimvar}) for $\gamma < \gamma_0$ is analogous.
This proves the convergence in \eqref{TypeIlim}.  The second statement of the theorem follows from
Proposition \ref{Vgamma}.  \hfill $\Box$

\medskip

Next, we consider  Gaussian RFs  with spectral density $f= f_{\rm II}$ in  \eqref{typeI/IIsp}:
\begin{equation} \label{typeIIsp}
f(x,y) = \frac{g(x,y)}{|x|^{2d_1} |y|^{2d_2}}, \qquad (x,y) \in \Pi^2.
\end{equation}

\begin{proposition}\label{farimaLM}
Let $X $
be a stationary Gaussian RF on $\Z^2$ with zero mean and spectral
density $f$ in  \eqref{typeIIsp}, where $0< d_1, d_2 < 1/2$ and  $g \ge 0 $ is a bounded function such that $\lim_{x, y  \to 0} g(x,y) =1$.
Then for any $\gamma >0 $
\begin{eqnarray}\label{convFBS}
n^{-H(\gamma)} S_{n\gamma} (x,y)
&\limfdd&\kappa(d_1)\kappa(d_2)
B_{d_1+.5,d_2+.5}(x,y), \quad (x,y) \in \R^2_+, \quad n \to \infty,
\end{eqnarray}
where $H(\gamma) := (1+ \gamma)/2 + d_1 + d_2 \gamma$ and
$B_{d_1 + .5, d_2 + .5}$ is a fractional Brownian sheet (see Example  \ref{exFBS} for definition), $\kappa^2(d) :=
\int_{\R} |1- \e^{\i x}|^2 |x|^{-2-2d} \d x = \pi (2 (d +  .5)^2 \Gamma(d) \cos(\pi d))^{-1}.$
As a consequence,  $X$ does not exhibit scaling transition
for any  $0< d_1, d_2 < 1/2 $.
\end{proposition}

\noi {\it Proof.} We follow the proof of Theorem \ref{thmGauss}.  Accordingly, it suffices to show \eqref{Hlimvar},
where
\begin{equation} \label{SpFBS}
V_\gamma(x,y) \ := \  \int_{\R^2}\frac{ (1-\e^{{\i} u x})(1 - \e^{{\i} v y})}{\i^2 u v}\, \sqrt{h(u,v)}  W(\d u, \d v),
\quad \text{with} \quad h(u,v) := |u|^{-2d_1} |v|^{-2d_2},
\end{equation}
is the spectral representation of fractional Brownian sheet, viz.,
$V_\gamma  \eqfdd \kappa(d_1) \kappa(d_2) B_{d_1+ .5, d_2+ .5}$, see  (\cite{leo2011}, (7))
and
$W(\d u, \d v)$ is the same as in \eqref{Vlimit00}.  Note the r.h.s. of \eqref{SpFBS} does not depend on $\gamma$.
Let $m= n^{\gamma}$. Then
\begin{eqnarray*}
&&R_{n \gamma}(x,y; x',y')  \\
&&=\ n^{-2H(\gamma)}\int_{\Pi^2} D_{[nx]}(u)\overline{ D_{[nx']}(u)} D_{[my]}(v)\overline{ D_{[my']}(v)}
|u|^{-2d_1}|v|^{-2d_2} g(u,v) \d u \d v \\
&&\sim\ n^{-1-2d_1}\int_\Pi D_{[nx]}(u)\overline {D_{[nx']}(u)}  |u|^{-2d_1}  \d u \
m^{-1-2d_2}\int_\Pi D_{[my]}(v)\overline{ D_{[my']}(v)}  |v|^{-2d_2}  \d v.
\end{eqnarray*}
The limit
$\lim_{n \to \infty} n^{-1-2d} \int_{\Pi} D_{[nx]}(u)\overline{ D_{[nx']}(u)} |u|^{-2d} \d u  =  (1/2)\kappa^2(d)
(x^{2d+1} + x'^{2d+1} - |x-x'|^{2d+1}), \ 0<d< 1/2,\,  x, x' >0 $
is  well-known.
This proves  \eqref{Hlimvar} and the proposition, too. \hfill $\Box$

\medskip

Given a line $\ell = \{ax+ by = 0\}\subset  \R^2$, a function $k(u,v), (u,v) \in \R^2_0 $ is said {\it  $\ell$-degenerated} if $k(u,v) =
\tilde k(au + bv), (u,v) \in \R^2_0$, where
$\tilde k(z)$ is a function of a single variable $z\in \R$. Obviously, $h(u,v)$ in  \eqref{h} is not
$\ell$-degenerated for any $\ell$. The proof of Lemma \ref{lemspec} is given in sec.4.

\begin{lemma} \label{lemspec} Let
\begin{equation} \label{Vfield}
V(x,y)
= \int_{\R^2}\frac{ (1-\e^{{\i} u x})(1 - \e^{{\i} v y})}{\i^2 u v}\, \sqrt{k(u,v)}   W(\d u, \d v), \quad (x,y) \in \bar \R^2_+,
\end{equation}
be a Gaussian RF, where $  W(\d u, \d v)$ is the same as in (\ref{Vlimit00}) and
$k(u,v)\ge 0$ is a measurable function such that
\begin{equation}\label{Kfin}
\int_{\R^2} \frac{k(u,v) \d u \d v} {(1+u^2)(1+v^2)}  \ < \  \infty.
\end{equation}
Let  $\ell$ be a  line in  $\R^2$.
Then $ V= \{V(x,y); (x,y) \in \bar \R^2_+\}$
has independent rectangular increments in direction $\ell$ if and only if $k$ is $\ell$-degenerated. Moreover,
$V$ does not have invariant rectangular increments in any direction.
\end{lemma}

\medskip

The following proposition obtains scaling transition of different type than Type I
for the class of stationary Gaussian RFs with
spectral density
\begin{equation}\label{Lavsp}
f(x,y) = g(x,y)|a x + b y|^{-2d},
\end{equation}
where $0< d < 1/2, a, b $ are parameters and $g$ is continuous at the origin.
Notice that when $ab =  0$,
\eqref{Lavsp} is the limiting case of Type II density in  \eqref{typeIIsp}   when
one of the parameters $d_1,  d_2 $ approaches zero. Partial sums limits of such RFs
where discussed in Lavancier \cite{lav2007}.

\begin{proposition}\label{Lav}
Let $X $
be a stationary Gaussian RF on $\Z^2$ with zero mean and spectral
density $f$ in  \eqref{Lavsp}, where $0< d < 1/2, ab \neq 0$ and  $g \ge 0 $ is a bounded function such that $\lim_{x, y  \to 0} g(x,y) =1$.
Then for any $\gamma >0 $ the limit $V_\gamma$ in \eqref{TypeIlim} exists and is written as in
\eqref{SpFBS} with  $H(\gamma)>0$ and $h(u,v) = h_\gamma(u,v)$ given by
\begin{equation}\label{hLav}
H(\gamma) := \begin{cases}(1+\gamma)/2 + d, &\gamma \ge 1,  \\
(1+ \gamma)/2 + \gamma d, &\gamma < 1,
\end{cases}
\qquad
h_\gamma(u,v) := \begin{cases}
|a u +  b v|^{-2d}, &\gamma = 1, \\
|a u|^{-2d}, &\gamma > 1, \\
|b v|^{-2d}, &\gamma < 1.
\end{cases}
\end{equation}
In particular, the RF $X$  exhibits scaling transition at $\gamma_0 =1 $ with $V_+ = \kappa_+ B_{1/2 + d,1/2}, \
V_- = \kappa_- B_{1/2, 1/2 + d}$, where $\kappa_\pm \ne 0$ are some constants.
Moreover, $X$ does not have Type I distributional
LRD property.

\end{proposition}

\noi {\it Proof.} Let $\gamma > 1$. Then similarly as in the proof of Theorem  \ref{thmGauss}, case
$\gamma > \gamma_0, 0< H_1 < 1$, with $m = n^\gamma $, we have that
$R_{n \gamma}(x,y; x',y')= \int_{\R^2} G_n(u,v) \d u \d v$, where $G_n$ is given in \eqref{Gnuv}  and
$n^{-2d} f(\frac{u}{n}, \frac{v}{m}) \to h_\gamma(u,v) = |au|^{-2d} $ a.e. in $\R^2$; moreover, $G_n$ is dominated by integrable function
$\bar G(u,v) := C(1+ u^2)^{-1} (1+ v^2)^{-1} |u|^{-2d}$. This proves
(\ref{Hlimvar}) for $\gamma > 1$ and the proof in the remaining cases $\gamma = 1$ and $\gamma < 1 $ is analogous.
The fact that $X$ admits  scaling transition at $\gamma_0 = 1$ is obvious since
$B_{1/2 + d,1/2} \neqfdd
c B_{1/2,1/2 + d} $ for any $c>0$. Finally,  since $h_1 $ in \eqref{hLav} is degenerated, by Lemma  \ref{lemspec}
the RF $V_{1} $ has independent increments
in the direction perpendicular to the line $a x + b y = 0$ and therefore $X$ does not have Type I distributional LRD property.
Proposition \ref{Lav} is proved. \hfill $\Box$

\begin{rem}\label{typeIII} {\rm  The above result can be described as an abrupt change of
the `dependence axis' of  RF $X$ under unbalanced
scaling. The form of spectral  density
in  \eqref{Lavsp} suggests that the LRD in $X$ is essentially `one-dimensional' along the
line $\ell = \{at + bs = 0\} \subset \R^2$. The `supercritical regime' $\gamma>1$  transforms this
`dependence axis' $\ell $ into the horizontal axis, since  $V_+ = \kappa_+ B_{1/2 + d,1/2}$ has
independent increments in the vertical direction. A similar transformation of $\ell$ into the
vertical axis occurs in the `subcritical regime' $\gamma < 1$.

}
\end{rem}

\begin{rem}\label{remON} {\rm  As noted in \cite{ps2014},  scaling transition occurs for a very different class of models
under joint temporal and contemporaneous
aggregation of {\it independent}
LRD processes in telecommunication and economics, see \cite{miko2002}, \cite{gaig2003}, \cite{domb2011},
\cite{pils2014} and the references therein. In these works, $\{X(t,s); t \in \Z\}, s \in \Z$ are independent copies of a stationary
LRD process $Y = \{Y(t); t \in \Z\}$
and the scaling limits $V_\gamma $ of the RF $X = \{X(t,s); (t,s)\in \Z^2\}$ necessarily have
independent increments in the vertical direction for any $\gamma >0$, meaning that $X$ cannot have
Type I distributional LRD by definition. Nevertheless for heavy-tailed centered ON/OFF process $Y$ and some other duration based
models, the results in   \cite{miko2002} imply that the above RF $X$ can exhibit  scaling transition with some $\gamma_0 \in (0,1)$ and
markedly distinct `supercritical' and `subcritical' unbalanced scaling limits $V_\pm$, viz., $V_+$ being a Gaussian RF with dependent increments in the horizontal direction
and $V_-$ having $\alpha-$stable $(1< \alpha < 2)$ distributions and independent increments in the  horizontal direction.
The well-balanced scaling limit $V_{\gamma_0} $ in the above models was discussed in detail in \cite{gaig2006}, \cite{pils2014} and was
shown to have interesting `intermediate' properties between $V_+$ and $V_-$.}
\end{rem}

\section{Proof of Lemma \ref{lemspec}}

We use
some facts about generalized functions (Schwartz distributions) (see \cite{yos1965}).
Let $S(\R^\nu) (\nu = 1,2)$ be the Schwartz space of all rapidly decreasing $C^\infty-$functions $\phi: \R^\nu \to  \R$, and
$S'(\R^\nu)$ be the space of all  generalized functions $T: S(\R^\nu) \to \R$.  The Fourier transform
$\widehat T \in S'(\R^2)$ of  $T \in S'(\R^2)$
is defined as $\widehat T(\phi) = T(\hat \phi),  \phi \in S(\R^2),$ where
$\widehat \phi (u,v) := \int_{\R^2} \e^{\i (ux + vy)} \phi(x,y) \d x \d y$. For $\varphi, \psi \in  S(\R)$,  denote
$(\varphi \otimes \psi)(u,v) := \varphi(u)\psi(v), \, (\varphi \otimes \psi) \in  S(\R^2)$.
Note that for any rectangle $K$ we have
$V(K) = \int_{\R^2}  \widehat \1_K(u,v) , \sqrt{h(u,v)}   W(\d u, \d v)$, where $\1_K $ is the indicator function of $K$.
It is easy to show that  $\{V(x,y)\}$ in (\ref{Vfield})
extends to a generalized stationary Gaussian random
field (\cite{dob1979}):
$$
{\cal  V}(\phi) := \int_{\R^2}\hat \phi(u,v) \sqrt{k(u,v)}  W(\d u, \d v), \qquad  \phi \in S(\R^2).
$$

Let $\ell $ be a given line and
$ \R^2_\pm (\ell')$ be the open halfplanes separated by line $\ell'  \perp \ell, 0\in \ell'$, viz.,
$\R^2_+ (\ell')\cup \ell' \cup  \R^2_- (\ell') = \R^2. $
Let us show
that the statements (a) $\{V(x,y); (x,y) \in \bar \R^2_+\}$ has independent rectangular increments in direction $\ell$ and (b) $  {\cal  V}(\phi_+)$ and  ${\cal  V}(\phi_-) $ are independent  for any
$\phi_\pm \in S(\R^2)$ with supports in  $\R^2_\pm (\ell')$
are equivalent. Statements (a) and (b) can be rewritten as
\begin{equation}\label{ortV1}
\int_{\R^2 } \hat \phi_+(u,v) \overline{\hat \phi_- (u,v)} k(u,v) \d u \d v \ = \ 0,
\qquad \phi_\pm  \in  L (\R^2_\pm(\ell')),
\end{equation}
and
\begin{equation}\label{ortV}
\int_{\R^2 } \hat \phi_+(u,v) \overline{\hat \phi_- (u,v)} k(u,v) \d u \d v \ = \ 0, \qquad \phi_\pm \in S(\R^2),
\qquad {\rm supp}(\phi_\pm)  \subset  \R^2_\pm (\ell'),
\end{equation}
respectively, where $L(\R^2_\pm (\ell')) $ is the set of all linear combinations of indicator functions
$\1_K $ of rectangles $K \subset \R^2_\pm (\ell')$.
To show the implication (b) $\Rightarrow$ (a),
note that any indicator functions $\1_{K_\pm},  K_\pm \subset \R^2_\pm (\ell')$
can be approximated in $L^2(\R^2) $ by
elements  $\phi_{\pm, \epsilon} \in S(\R^2), \epsilon >0$ with
compact supports ${\rm supp}(\phi_{\pm,\epsilon})
\subset  \R^2_\pm (\ell')$.
The approximating functions can be taken as $\phi_{\pm,\epsilon} = \1_{K_\pm}
\star \theta_\epsilon, \epsilon >0,$ where
$\theta_\epsilon (u,v) := \epsilon^{-2} \theta(u/\epsilon, v/\epsilon)$, $\theta $ is a $C_0^\infty(\R^2)$ probability kernel, and $\star $ denotes the convolution. See  \cite{ste1971}, Thm.1.18.  Using
$|\widehat \phi_{\pm,\epsilon}(u,v)| = |\widehat \1_{K_\pm}(u,v)|\, |\widehat \theta_\epsilon(u,v)|, \, |\widehat \1_{K_\pm}(u,v)| \le
C(1+ |u|)^{-1}(1+|v|)^{-1}$ and $|\widehat \theta_\epsilon(u,v)| = |\widehat \theta (\epsilon u,
\epsilon v)| \le C $, with $C<\infty$ independent of $\epsilon \to 0$, we can easily show that
\eqref{ortV} implies \eqref{ortV1} for $\phi_\pm = \1_{K_\pm}$ and any rectangles
$K_\pm \subset \R^2_\pm (\ell')$, proving the implication  (b) $\Rightarrow$ (a).

Next, consider the converse implication (a) $\Rightarrow$ (b). It suffices to prove
\eqref{ortV} for $\phi_\pm \in S(\R^2)$ with  compact supports
${\rm supp}(\phi_\pm)  \subset  \R^2_\pm (\ell')$. Consider approximation of such $\phi_\pm$
by  step functions
$\phi_{\pm,n}(x,y) := \sum_{(k,j) \in \Z^2} \phi(k/n,j/n) $  $\1_{(k/n, (k+1)/n] \times (j/n, (j+1)/n]} (x,y), \,
n \in \N_+$.
Then $\sup_{x,y} |\phi_{\pm,n}(x,y) - \phi_{\pm}(x,y)| \to 0 \
(n \to \infty)$ and
$\phi_{\pm,n} \in L(\R^2_\pm (\ell')) $ for all $n > n_0 $ large enough;
moreover,
\begin{equation}\label{phibdd}
|\widehat \phi_{\pm,n}(u,v)| \le C(1+ |u|)^{-1}(1+|v|)^{-1}
\end{equation}
with $C<\infty $ independent of $n$. Inequality \eqref{phibdd} can be shown
using summation by parts, as follows. We have
$|\widehat \phi_{\pm,n}(u,v)| = $   $|uv|^{-1}\big|\sum_{k,j} \psi_{\pm,n}(k,j)
\e^{\i (ku/n)} \e^{\i(jv/n)} |
=  |uv|^{-1}\big|\sum_{k,j} \psi_{\pm, n}(k,j) (\e^{\i (ku/n)}-1) (\e^{\i(jv/n)}-1) |$, where
$\psi_{\pm,n}(k,j) := \phi_\pm ((k-1)/n, (j-1)/n) - \phi_\pm ((k-1)/n, j/n) -
\phi_\pm (k/n, (j-1)/n) + \phi_\pm (k/n, j/n)$ is the double difference satisfying
$ |\psi_{\pm, n}(k,j)| \le C/n^2 $ (recall that $\phi_\pm $ is infinitely differentiable with compact support).
Hence,  
$|\widehat \phi_{\pm,n}(u,v)|
\le C|uv|^{-1} (1 \wedge |u|) (1\wedge |v|) \le C (1+ |u|)^{-1}(1+|v|)^{-1}$,
proving \eqref{phibdd}. The above facts together with condition \eqref{Kfin}
allow  using  the dominated convergence criterion to prove \eqref{ortV}
from \eqref{ortV1} via the above approximation by step functions.

Consider first the case when $\ell$ is the horizontal axis,
in which case $\ell'$ is the vertical axis.  Let $k$ be $\ell$-degenerated,
or $k(u,v) = \tilde k(v) $ for some measurable function $\tilde k$ satisfying \eqref{Kfin}.
Let $\phi_\pm = \varphi_\pm \otimes \psi_\pm$, where  $\varphi_\pm, \psi_\pm \in S(\R)$
satisfy ${\rm supp}(\varphi_+) \subset (0, \infty), {\rm supp}(\varphi_-) \subset (-\infty, 0).$ Then $\E [{\cal  V}(\phi_+){\cal  V}(\phi_-) ]
= J_1 J_2 = 0$ follows from
$$
J_1 \ := \
\int_{\R} \hat \psi_+(v)
\overline{\hat \psi_- (v)}\tilde  k(v)\d v, \qquad J_2 :=
\int_{\R} \hat \phi_+(u)
\overline{\hat \phi_- (u)} \d u =  2\pi \int_{\R} \phi_+(z) \phi_- (z) \d z  = 0
$$
by Parseval's identity. Then \eqref{ortV} follows by taking
linear combinations of the above  $\phi_\pm$  in a standard way,  proving  that
$\{V(x,y);  (x,y) \in \bar \R^2_+\}$ has independent rectangular increments in the horizontal $\ell$.

Let us prove  the converse implication, i.e. that \eqref{ortV} implies that
$k$ is $\ell$-degenerated.  Since the bilinear form
${\cal T} (\hat \phi_+, \hat \phi_-) :=  \int_{\R^2 } \hat \phi_+(u,v) \overline{\hat \phi_- (u,v)} k(u,v) \d u \d v$
in (\ref{ortV}) is invariant with respect to shifts of $\phi_\pm $ in the horizontal direction, it is easy to see that
(\ref{ortV}) holds for any $\phi_\pm \in S(\R^2)$ with
\begin{equation}\label{supports}
{\rm supp}(\phi_+) \subset (\R\setminus [-\epsilon, \epsilon])\times \R, \qquad
{\rm supp}(\phi_-) \subset (-\epsilon, \epsilon)\times \R, \qquad  \epsilon>0.
\end{equation}
Let $\phi_+ = \varphi \otimes \psi_+, \phi_- = \phi_\epsilon \otimes \psi_-,$ where $\varphi, \phi_\epsilon, \psi_\pm \in S(\R)$
satisfy ${\rm supp}(\varphi) \subset \R \setminus \{0\}$ and $\phi_\epsilon(x) = \epsilon^{-1} \phi(x/\epsilon)$, where
$\phi \in S(\R)$ is a symmetric probability density with  ${\rm supp}(\phi) \subset (-1,1)$.  Note $\phi_\pm$ satisfy (\ref{supports}) for
$\epsilon>0$ small enough, implying
\begin{equation} \label{Tzero}
T_\epsilon(\hat \varphi) \ := \
{\cal T} (\hat \varphi \otimes  \hat \psi_+,  \hat \phi_\epsilon \otimes \hat \psi_-) \ = \  0,
\qquad {\rm supp}(\varphi) \subset \R \setminus \{0\}
\end{equation}
according to  (\ref{ortV}). We claim that $\lim_{\epsilon \to 0} T_\epsilon(\hat \varphi) = T(\hat \varphi)$, where
the generalized function $T \in  S'(\R)$ is given by
\begin{equation} \label{Tfunk}
T(\varphi) :=  \int_{\R} \varphi(u) \mathcal{K}(u)  \d u,  \quad
  \quad \mathcal{K}(u) :=
\int_{\R} \hat \psi_+ (v) \overline{\hat \psi_- (v)} k(u,v) \d v.
\end{equation}
Indeed, $|T_\epsilon(\hat \varphi)- T(\hat \varphi)|  \le \int_{\R^2} |\hat \varphi(u)| |\hat  \phi_\epsilon (u) - 1|
|\hat \psi_+(v)| |\hat \psi_-(v)|  |k(u,v)| \d u \d v  \to 0 $ when   $\epsilon \to 0$,
which follows from  $\hat  \phi_\epsilon (u) \to 1$ and the integrability of $ \hat \varphi(u)
\hat \psi_+(v) \hat \psi_-(v) k(u,v) $ on $\R^2$; see \eqref{Kfin},
proving the above claim.
Hence and from \eqref{Tzero}, we infer that
$\hat T(\varphi) = T(\hat \varphi) $ vanishes for $\varphi$ with ${\rm supp}(\varphi) \subset \R \setminus \{0\}$.
The last fact implies that $\hat T $
is a linear combination
of Dirac's  $\delta$-function and its derivatives:  $\hat T(\varphi) = \sum_{k=0}^p c_k \varphi^{(k)}(0), \, \varphi^{(k)}(x) := \d^k \varphi(x)/\d x^k$
with real coefficients $c_k, 0\le k \le p <\infty,$ see
(\cite{yos1965}, Ch.1.13, Thm.3). In other words,
$T(\varphi) = (2\pi)^{-1} \sum_{k=0}^p c_k \hat \varphi^{(k)}(0)
=  \int_{\R} P(u) \varphi (u) \d u, $ where $P(u) =  (2\pi)^{-1} \sum_{k=0}^p c_k u^k $ is a polynomial
of degree $p$.  Comparing the last expression with (\ref{Tfunk}) we obtain that $ \mathcal{K}(u) = P(u) $ a.e. in $\R$. Since
$|\mathcal{K}(u)| \le C \int_{\R} (1+v^2) k(u,v) \d v  $ satisfies $\int_{\R}  |\mathcal{K}(u)| (1+u^2)^{-1} \d u < \infty $, see
 \eqref{Kfin}, this means that the function $\mathcal{K}$ in (\ref{Tfunk})  is constant on the real line: $\mathcal{K}(u) = (2\pi)^{-1} c_0. $
 Since the last fact holds for arbitrary  $\psi_\pm \in S(\R)$, we conclude that the function
 $k$ does not depend on $u $, viz., $k(u,v) = \tilde k(v)$ a.e. in $\R^2$. This proves the first statement of the lemma
for a horizontal line $\ell $.

The case of a general line $\ell = \{au+ bv = 0\}$ can be reduced to that of the horizontal line
$\ell_1 := \{ v = 0 \}$
by a rotation of the plane.  Indeed,
in such a case, similarly as above we can show that there exists an orthogonal $2\times 2-$matrix $O $
mapping $\ell $ to $\ell_1$
such that
$T_O(\hat \varphi) :=  \int_{\R} \hat \varphi(u)  {\mathcal K}_O(u)  \d u $ vanishes
for $\varphi$ with ${\rm supp}(\varphi) \subset \R \setminus \{0\}$,
where ${\mathcal K}_O$ is defined
as in (\ref{Tfunk}) with $k(u,v) $ replaced by $k_O(u,v) := k(O^{-1}(u,v))$, and consequently $k_O(u,v)$
is $\ell_1$-degenerated,  or $k(u,v)$ is $\ell$-degenerated.

Let us prove the second statement of the lemma. Assume {\it ad absurdum} that
$\{V(x,y); (x,y) \in \bar \R^2_+\}$ has invariant rectangular increments in the horizontal direction.
Then ${\cal T} (\hat \theta_a \phi, \hat \phi) = {\cal T} (\hat \theta_0 \phi, \hat \phi) $
does not depend on $a \in \R$, where $\theta_a \phi (x,y) := \phi (x+a, y) $ is a shifted function.  Since
${\cal T} (\hat \theta_a \phi, \hat \phi) =   \int_{\R^2 } \e^{\i a u}  |\hat \phi(u,v)|^2 k(u,v) \d u \d v \ \to \ 0 \, (a \to \infty) $ by the
Lebesgue theorem, we obtain a contradiction. The case when  $\{V(x,y);  (x,y) \in \R^2_+\}$ has invariant
rectangular increments in  arbitrary  direction can be
treated analogously.
Lemma \ref{lemspec}
is proved.

\section*{Acknowledgements}

This research was supported by a grant (no.\, MIP-063/2013) from the Research Council of Lithuania.


\end{document}